\documentclass{amsart}

\usepackage{amssymb}

\usepackage{graphicx}

\usepackage{blindtext}
\usepackage{enumitem}
\usepackage[colorlinks]{hyperref}
\usepackage{etoolbox,lastpage}
\usepackage{tikz}
\usetikzlibrary{arrows}

\newtheorem{theorem}{Theorem}[section]
\newtheorem{lemma}[theorem]{Lemma}
\newtheorem{corollary}[theorem]{Corollary}

\theoremstyle{definition}
\newtheorem{definition}[theorem]{Definition}
\newtheorem{example}[theorem]{Example}

\theoremstyle{remark}
\newtheorem{remark}[theorem]{Remark}

\numberwithin{equation}{section}

\begin{document}

\title[A new method to prove the collatz conjecture]{A new method to prove the collatz conjecture}
\author[D.Karami]{Danial Karami}
\address{}
\email{}
\thanks{}

%    \subjclass is required.

\date{}

\dedicatory{}


%    Abstract is required.
\begin{abstract} The Collatz conjecture is a famous math problem that was introduced by Lothar Collatz in 1937, and nobody has yet succeeded in proving or disproving it. In this article, I will analyze this problem with a new approach and I will discuss my important findings that are helpful in getting closer to the proof of this unsolved problem in mathematics.
\end{abstract}

\maketitle

\section{Introduction} The Collatz conjecture, also know as Syracus conjecture, involves a sequence where we start with any positive integer. If it is odd, triple it and add one to get the next number in the sequence. However, if it is even, divide it by 2 instead. The conjecture claims that the value of the original selected number is not important because the sequence will always reach one.

As of today, this conjecture has been checked by computer tests up to $2.95\times10^{20}$ and most of the mathematicians believe that this conjecture is true for all natural numbers. However, nobody has yet succeeded in providing a proof that can convince us of whether this is correct or not.

\begin{lemma} The configuration of the Collatz conjecture can be viewed as

\centerline{any chosen number $\xrightarrow{\text{\hspace{10 mm}}}$1}
\end{lemma}

\begin{theorem} We can define this conjecture as a function f:

${\displaystyle f(n)={\begin{cases}n/2&{\text{if }}n\equiv 0{\pmod {2}}\\3n+1&{\text{if }}n\equiv 1{\pmod {2}}\end{cases}}}$

In fact, we form a sequence by performing this action many times. Begin with any natural number and take the result of each step as the input for the next step, In notation:

 ${\displaystyle a_{i}={\begin{cases}n&{\text{for }}i=0\\f(a_{i-1})&{\text{for }}i>0\end{cases}}}$

Where ${\displaystyle a_{i}}$ is equivalent to ${\displaystyle f}$ being applied to ${\displaystyle n}$ ${\displaystyle i}$ times (${\displaystyle a_{i}=f^{i}(n)}$)~\cite{1}.
\end{theorem}

\begin{example} Here is an example with 17 as the selected first number:

$17\rightarrow52\rightarrow26\rightarrow13\rightarrow40\rightarrow20\rightarrow10\rightarrow5\rightarrow16\rightarrow8\rightarrow4
\rightarrow2\rightarrow1$
\end{example}

\begin{definition} \label{1.4}
Below, the modified form of Collatz function is discussed:

${\displaystyle f(n)={\begin{cases} n/2&{\text{if }}n\equiv 0\\(3n+1)/2&{\text{if }}n\equiv 1\end{cases}}{\pmod {2}}}$ ~\cite{2}

It is clear that $3n+1$ is always an even number when $n$ is odd such that $n$ $\in$ $\mathbb{N}$, and this is the reason why we can write it in this way.
\end{definition}

\textbf{2-adic extension} "Consider the following function $T$:

${\displaystyle T(x)={\begin{cases}x/2&{\text{if }}x\equiv 0{\pmod {2}}\\(3x+1)/2&{\text{if }}x\equiv 1{\pmod {2}}\end{cases}}}$

It is clear that due to the 2-adic measure, the loop of 2-adic integers, $\mathbb {Z} _{2}$, is persistent and preserves measure. In fact, its dynamics has recognised to be ergodic. We can show the parity vector function ${\displaystyle Q}$ which perform on $\mathbb {Z} _{2}$ as:

${\displaystyle Q(x)=\sum _{k=0}^{\infty }\left(T^{k}(x)\mod 2\right)2^{k}}$

 As a result, the function $Q$ is a 2-adic isometry. All of the infinite and unlimited parity sequence occur for one 2-adic integer so that in $\mathbb {Z} _{2}$, almost all of the routes are acyclic. An equivalent formulation of this problem in mathematics is ${\displaystyle Q(\mathbb {Z} ^{+})\subset {\frac {1}{3}}\mathbb {Z} }$." This area has been widely studied in~\cite{3}.

\textbf{Iteration on the real or complex digits}: "As restriction to the numbers of the smooth real and complex map, Collatz conjecture map can be seen as

${\displaystyle f (z) = {\frac {1}{2}}z.\cos^{2}\left({\frac {\pi }{2}}z\right) + (3z+1).\sin ^{2}\left({\frac {\pi }{2}}z\right)}$

If it is the standard Collatz map explained above, we can improve it by replacing $3n+1$ with $(3n+1)/2$. Therefore, this map can be shown as

${\displaystyle f (z) = {\frac {1}{2}}z.\cos^{2}\left({\frac {\pi }{2}}z\right) + {\frac {(3z+1)}{2}}.\sin ^{2}\left({\frac {\pi }{2}}z\right)}$

Chamberland examined the repetition on the real number line in 1996. He demonstrated that when there are infinitely many fixed points, the Collatz conjecture doesn't hold for the real numbers. Furthermore, the orbits and cycles escape monotonously and are boundless." For more details study~\cite{4}.

\begin{definition}
 "The result of $3k+1$ is always equal to an even number, unless the case is $k$ = 2\emph{\H{k}}. So with $a$ $\geq$ 1 and \emph{\'{k}} odd, $3k+1$ = $2^{a}$\'{k}. From the set $\textsl{I}$ of odd integers into function f such that $f(k) = \acute{k}$, syracuse function will appear. Some properties of this function are:

\item $\bullet$  $\forall$ $k$ $\in$ $I$ ; $f(4k+1)$ = $f(k)$. Because $3(4k+1)+1 = 12k+4 = 4(3k+1)$

\item $\bullet$ We can generalize more: $\forall$ $p$ $\geq$ 1, $h = 2k+1$ ; $f^{p-1}$($2^{p}$$h-1$) = 2 $\times$ $3^{p-1}$h-1.
(Here $f^{p-1}$ is function iteration notation)

\item $\bullet$ $\forall$ $h = 2k+1$ ; $f(2h-1)$ $\leq$ {$\frac{3h-1}{2}$}

The Collatz conjecture is commensurate to express that, for every $k$ $\in$ \textit{I}, there is an integer $n$ $\geq$ 1, s.t. $f^{n}$(k) = 1." More information can be found in ~\cite{5}
\end{definition}

 Anyway, in the main result of this paper, particular patterns of the numbers that approach 1 are mentioned. In fact, it has been discovered that in contrast with the fact that each number has its own algorithm, they have to follow some same particular edict in their journey to 1.

  \section{Main result}
In order to find the steps before $1$, lets lead 1 to other numbers by doing the formulation of the conjecture backward.
\begin{theorem}\label{2.1}
One step before $1$, there are numbers of the form $2^{n}$, where $n$ $\in$ $\mathbb{N}$.
\begin{proof}
If we reach a number of the form $2^{n}$, we will keep dividing it by 2 until we reach 1. It is also guaranteed that we will reach
$2^{1}$ = 2 as it is the step before 1. Note that we are not able to reach 1 with 3n + 1 because then n = 0, which should never occur! In summary:
$2^{n}$ $\xrightarrow{\text{$2^{n}\div2$}}$ $2^{n-1}$ $\xrightarrow{\text{$2^{n-1-1-...}$(for n times)$$}}$ $2^{n-n}$ = $2^{0}$ = 1
\end{proof}
\end{theorem}
\begin{definition} \label{2.2}
Let "$2^{n}$s" be numbers of the form $2^{n}$. In this conjecture, the numbers before 1 are in two categories. The first category is \{$2^{n}$ $\mid$ n $\in$ $\mathbb{N}$\}, while the second one is \{$x$ $\in$ $\mathbb{N}$ $\mid$ $x$$\neq$$2^{n}$, $x$$\neq$$1$\}.

The numbers in the first group reach 1 just after being divided by 2 for $n$ times (proof in theorem~\ref{2.1}). As a result, it is necessary to find a connection between the numbers in the second category with $2^{n}$s since $2^{n}$s arrive at 1 directly.

Now let us try leading $2^{n}$s to other numbers by doing the formulation of the conjecture backward. Multiplying a $2^{n}$ by 2 for $i$ times only yields another number of the form $2^{n}$ and this is useless. Hence, the only way to relate $2^{n}$s to other numbers is to subtract 1 from $2^{n}$s, and then divide them by 3 (except $2^{2}$ because $\frac{4-1}{3}$ = 1 and 1 is our destination, not a step before $2^{n}$s).
\end{definition}

\textbf{Proposition} To identify the step before $2^{n}$s, we can divide any $2^{n}-1$ by 3. Always the result would be an odd number.
\begin{proof} We proceed with proof by contradiction. Thus, instead of considering that $(2k+1)/3$ = 2\'{k} + 1, lets consider that $(2k+1)/3$ $\neq$ 2\'{k} + 1 so $(2k+1)/3$ $=$ \emph{2\H{k}}. Now we have
$\frac{2k+1}{3}$ = \emph{2\H{k}} $\Rightarrow$ $2k+1$ $=$ \emph{6\H{k}} $\Rightarrow$ $2k+1$ = 2(\emph{3\H{k}}) $\Rightarrow$ $2k+1$ = $2k^{\prime \prime \prime}$

Exactly the last stage has a contradiction. It shows that $\neg$($\frac{2k+1}{3}$= 2\'{k}+1) is false because 2k+1 $\neq$ $2k^{\prime\prime\prime}$, so the phrase $\frac{2k+1}{3}$ = 2\'{k}+1 is true.
\end{proof}

\begin{remark} \label{2.3}
There is an important question regarding definition~\ref{2.2}. Is every number of the form $2^{n}-1$ divisible by 3? To answer this question, let’s first analyze the following algorithm for having a better understanding:

\{$2^{1}-1$ $\neq$ 3$a$, $2^{2}-1$ = 3$b$, $2^{3}-1$ $\neq$ 3$c$, $2^{4}-1$ = 3$d$, $2^{5}-1$ $\neq$ 3$e$, ...\}

As it can be seen in the above set of numbers, just $2^{2k}-1$, which is $4^{m}-1$, is divisible by 3, not $2^{2k+1}-1$. Lets confirm this claim!
\end{remark}

\textbf{Proposition} $2^{2k}-1$ = $4^{m}-1$ is divisible by 3, for any integer $n$ $\in$ $\mathbb{N}$.

\begin{proof}
Let $S_m$ = $4^{m}-1$ = 3$r$ where $r$ $\in$ $\mathbb{N}$. In the condition that $k$ = $m$, we would have $S_k$ = $4^{k}-1$ = 3$r$. We must show that for the next turn, which is $S_{k+1}$, $4^{k+1}-1$ = 3$t$ such that $t$ $\in$ $\mathbb{N}$. The reason is that $m$ represents the natural numbers.

Lets multiply $S_k$ by 4 and then add 3 to the consequence:

  4 $\times$ ($4^{k}-1)+3$ = 4 $\times$ $(3r)+3$ $\Rightarrow$ $4^{k+1}$-4+3 = $12r+3$ $\rightarrow$  $4^{k+1}$-1 = $12r+3$

By factoring out a 3, the result can be viewed as

  $4^{k+1}-1$ = $3(4r+1)$ $\Rightarrow$ $t$ = $4r+1$

Clearly, $4^{k+1}-1$ = $3(4r+1)$ is divisible by 3, and the proof is over.
\end{proof}

\textbf{Proposition} $2^{2k+1}-1$ is not divisible by 3 for any integer $k$ $\in$ $\mathbb{N}$.

\begin{proof} According to the previous proof, $2^{2k}-1$ = $3r$ and it is divisible by 3 such that $r$ $\in$ $\mathbb{N}$. Anyway lets multiply the whole equality by 2:

2 $\times$ $(2^{2k}-1)$ = 2 $\times$ $(3r)$ $\Rightarrow$ $2^{2k+1}-2$ = $6r$

Then, we can add 1 to the both sides of the equality.

$2^{2k+1}-2+1$ = $6r+1$ $\Rightarrow$ $2^{2k+1}-1$ = $6r+1$

As a result, $2^{2k+1}$-1 is not divisible by 3 because $\frac{6r+1}{3}$ $\neq$ $j$ such that $j$ $\in$ $\mathbb{N}$. Here the proof comes to an end.
\end{proof}

Based on the details written in remark~\ref{2.3}, the odd number $2k+1$ is the step before $2^{n}$s if an only if 3 $\times$ $(2k+1)+1$ = $2^{2k}$ = $4^{m}$. However, there is an exception which is 1, and the reason is that we do not perform any math operation on it since 1 is our final target in this conjecture.

\begin{definition} \label{2.4}
According to all of the above discussions, the simplified template of the Collatz conjecture can viewed as

... $\rightarrow$ \emph{2\'{k}} $+ 1$ $ \xrightarrow{\text{$\times3+1$}}4^{m} = 2^{n} $$\xrightarrow{\text{($\div2$) for n time}}1$

Furthermore, it should be pointed out that odd numbers are the result of dividing the even numbers by 2 until it is not possible. Because if $3n+1$ = \emph{2\H{k}}$ + 1$, $n$ = $2t$ and this is not acceptable due to the rules of the Collatz conjecture which prevent us to do $3n+1$ when $n$ is even. All in all, it can be concluded that there is an even number before any odd number in this conjecture; so lets improve the above configuration by adding a $2k$ before \emph{2\'{k}} $+ 1$:

... $\rightarrow$ $2k\xrightarrow{\text{($\div 2$) for several times}}$ \emph{2\'{k}} $+ 1$ $ \xrightarrow{\text{$\times3+1$}}4^{m} = 2^{n} $$\xrightarrow{\text{($\div2$) for n time}}1$

Here we mentioned that there always exist a $2k$ before any \emph{2\'{k}}$ + 1$. So currently we need to be able to prove why every even number will always reach an odd number after being divided by 2 for many times.
\end{definition}

\begin{theorem}\label{2.5} All even numbers will eventually arrive at odd numbers.

\begin{proof} To begin, choose a natural even number and divide it by 2 once. Assuming that $2k$ $\div$ 2 = \emph{2\'{k}}$ + 1$, we are done because the result is an odd number fortunately. But, if $2k$ $\div$ 2 = \emph{2\H{k}}, we have to keep dividing it by 2 until we reach an odd number. In fact, even numbers can be defined as $2k$ = $2^{n}$ $\times$ $x$ where $x$ = \emph{2\'{k}} $+ 1$. Here $n$ describes how many times the $2k$ needs to be divided by 2 until it reaches an odd number like $x$. For instance, 48 will reach 3, after we divide it 4 times by 2; in notation: 48 = $2^{4}$ $\times$ $3$
\end{proof}

\end{theorem}

\begin{remark} \label{2.6}
 Pursuant to the content of definition ~\ref{2.4}, up to 2 stages before $2^{n}$s ($4^{m}$s) are going to be specified in this section.

 First of all, find all values of $\frac{4^{m}-1}{3}$, the reverse form of $3n+1$, for all $m$ $\in$ $\mathbb{N}$ albeit we must draw our attention to exclude $\frac{4^{1}-1}{3}$ since 1 is our destination. A set of numbers can formed which can be represented as

$\mathbb{A}$ = \{5, 21, 85, 341, 1365, 5461, ... \}.

Whenever $2k+1$ $\xrightarrow{\text{$\times3+1$}}4^{n}$, odd number $2k+1$ is surely a member of set \textbf{A}.

 We can further examine the second stage before $2^{n}$s. As we discussed in definition~\ref{2.4}, there exists an even number before any odd number in this conjecture. So by multiplying numbers in \textbf{A} by 2 (backward form of the n/2), we can form a sequence, which is a step before the set of \textbf{A}, and sequentially define it as:

$\mathbb{B}$ = \{10, 42, 170, 682, 2730, 10922, ... \}

So if \emph{2\'{k}} $\xrightarrow{\text{$\div2$}}$ $2k+1$ where $2k+1$ $\in$ \textbf{A}, even number \emph{2\'{k}} is a member of set \textbf{B}.

\end{remark}

\begin{theorem} \label{2.7}
Sets \textbf{A} and \textbf{B} follow their own rule to get the next number. The algorithm of set \textbf{A} can be defined as

$\mathbb{A}$ = $\{$5$\xrightarrow{\text{$\times4+1$}}21$$\xrightarrow{\text{$\times4+1$}}85$
$\xrightarrow{\text{$\times4+1$}}341$$\xrightarrow{\text{$\times4+1$}}...$\}$\Rightarrow$ $t_n$ = $\frac{4^{n+1}-1}{3}$

And the algorithm of set \textbf{B} is defined as

$\mathbb{B}$ = $\{$10$\xrightarrow{\text{$\times4+2$}}42$$\xrightarrow{\text{$\times4+2$}}170$
$\xrightarrow{\text{$\times4+2$}}682$$\xrightarrow{\text{$\times4+2$}}...$\}$\Rightarrow$ $t_n$ = $\frac{2(4^{n+1}-1)}{3}$

\end{theorem}

\begin{corollary} \label{2.8}
  From the results obtained in this paper until now, the configuration of this problem can be illustrated as

 ... $\rightarrow$ \{$x$ $\mid$ $x$ $\in$ $\mathbb{B}$\} $\rightarrow$ \{$y$ $\mid$ $y$ $\in$ $\mathbb{A}$\} $\rightarrow$ $4^{n}$ $\rightarrow$ 1

 However this is a general style of this conjecture, so it is intellectual that probably your entry number could be in middle of the above configuration.
\end{corollary}

Anyhow, lets clarify and summarize the results of this paper in the below section so that we can have a better understanding of the Collatz conjecture!

\begin{corollary} \label{2.9}

The starting integer can be an even number or an odd number. Assuming that the starting integer is even, it will eventually reach an odd number according to theorem~\ref{2.5}. Hence, it is better to consider that the consequence of primary steps is certainly an odd number. If the obtained odd number is 1, we have just reached our destination. On the other hand, if the acheived odd number is in \textbf{A}, it will arrive at $4^{m}$ = $2^{n}$ in the next step (proof in remark~\ref{2.6}) and then it will arrive at 1 easily as we proved in theorem~\ref{2.1}. In condition that neither is the case, by tripling it and adding 1, we will gain an even number which will lead us to another odd number with respect to theorem~\ref{2.5}.

All in all, this loop between odd and even numbers before the sets \textbf{A} or \textbf{B}, will repeat again and again. It is not going to stop until it reaches one of the numbers in \textbf{A} or \textbf{B}. So by considering all the results, the general and final form of the Collatz conjecture can be viewed as

\item \begin{tikzpicture}[->,>=stealth',auto,node distance=2.3cm,
  thick,main node/.style={circle,draw,font=\sffamily\small\bfseries}]

  \node (1) {2k};
  \node (2) [right of=1] {2\'{k}+1};
  \node (3) [right of=2] {\{$x$ $\mid$ $x$$\in$$\mathbb{B}$\}};
  \node (4) [right of=3] {\{$y$ $\mid$ $y$$\in$$\mathbb{A}$\}};
  \node (5) [right of=4] {$4^{m}$ = $2^{n}$};
  \node (6) [right of=5] {1};
  \path[every node/.style={font=\sffamily\small}]

    (1) edge[bend right] node [right] {} (2)
    (2) edge[bend right] node [left] {} (1)
    (2) edge node [right] {} (3)
    (3) edge node [right] {} (4)
    (4) edge node [right] {} (5)
    (5) edge node [right] {} (6);
\end{tikzpicture}

In case that the starting integer is not an even number, omit 2k from the above configuration. However, the loop between odd and even number is undeniable, so by making a little change in the above configuration, we get

\item \begin{tikzpicture}[->,>=stealth',auto,node distance=2.3cm,
  thick,main node/.style={circle,draw,font=\sffamily\small\bfseries}]

  \node (1) {2k+1};
  \node (2) [right of=1] {2\'{k}};
  \node (3) [right of=2] {\{$x$ $\mid$ $x$$\in$$\mathbb{B}$\}};
  \node (4) [right of=3] {\{$y$ $\mid$ $y$$\in$$\mathbb{A}$\}};
  \node (5) [right of=4] {$4^{m}$ = $2^{n}$};
  \node (6) [right of=5] {1};
  \path[every node/.style={font=\sffamily\small}]

    (1) edge[bend right] node [right] {} (2)
    (2) edge[bend right] node [left] {} (1)
    (2) edge node [right] {} (3)
    (3) edge node [right] {} (4)
    (4) edge node [right] {} (5)
    (5) edge node [right] {} (6);
\end{tikzpicture}

Generally, the above configurations are telling us that as time goes by in the calculations using the Collatz function, numbers that exist in the loop between odd numbers $2k+1$ and even numbers \emph{2\'{k}}, gradually become smaller until they reach 1. However, sometimes numbers get bigger intermittently as they approach 1.
\end{corollary}
\begin{corollary} Because each number has its own path to 1, we can not write any specific algorithm that is appropriate for all numbers. But, it is logical that alternation between odd and even numbers is providing a good chance for the selected integer to reach one of the special class of numbers that we mentioned in this paper, and after that, everything is clear about the way it takes to arrive at 1.
\end{corollary}

 From the results of this article, we conclude that because the result of the primary steps is odd, if somebody were to explain that how odd integers finally arrive at numbers in \textbf{A} or \textbf{B}, the conjecture will be proven. In fact, the thing which need to be proven is that, how the loop between odd and even numbers introduced in corollary~\ref{2.9}, would reach the numbers in set \textbf{A} or \textbf{B}. So far, this has been the last approach to prove the Collatz conjecture.

\subsection*{\textbf{Acknowledgements}}I am really grateful to all of those who taught me to learn and i was honored to meet them.

\section*{\textbf{References}}

\begin{enumerate}
  \bibitem{1} J.J. O'Connor and E.F. Robertson, Lothar Collatz, St Andrews University School of Mathematics and Statistics, Scotland; (2006).
  \bibitem{2} Livio Colussi, The convergence classes of Collatz function, Theoretical Computer Science; vol.412 issue 39 pp. 5409–5419 (2011).
  \bibitem{3} Jeffrey C. Lagarias and Daniel J. Bernstein, The 3x+1 conjugacy map, Canadian journal of mathematics; vol.48 issue 6 pp. 1154-116 (1996).
  \bibitem{4} Marc Chamberland, A continuous extension of the 3x+1 problem to the real line, dynamic of continuous, discrete impulse system; vol.2 no.4 pp. 495-509 (1996).
  \bibitem{5} Jeffery C. Lagaries, The 3x+1 problem and its generalization, Journal of American Mathematical Monthly; vol.92 no.1 pp.3-23 (1985)
  \item Shalom Eliahou, The 3x+1 problem: new lower bounds on nontrivial cycle lengths, Discrete Mathematics; vol.118 issues 1-3 pp. 45–56 (1993).
  \item Leavens, Gary T., Vermeulen, Mike, 3x+1 Search Programs, Computers and Mathematics with Applications; vol.24 issue 11 pp. 79–99 (1992).
  \item Ben-Amram, Amir M., Mortality of iterated piecewise affine functions over the integers: decidability and complexity, Computability; vol.4 no.1 pp. 19-56 (2015), doi:10.3233/com-150032.
  \item Alex V. Kontorovich, Yakov G.Sinai, Structure theorem for (d,g,h)-Maps, https://arXiv.org/abs/math/0601622; (2006).
  \item G. H. Hardy and E. M. Wright, An Introduction to the Theory of Numbers, Oxford university; (2008).
\end{enumerate}
\end{document}